\documentclass{amsart}

\advance \topmargin by -\headheight
\advance \topmargin by -\headsep

\usepackage{graphics}
\usepackage{amscd}
\usepackage{amssymb}
\usepackage{longtable}

\newtheorem{thm}{Theorem}[section]
\newtheorem{prop}[thm]{Proposition}
\newtheorem{lem}[thm]{Lemma}
\newtheorem{cor}[thm]{Corollary}

\newtheorem*{thm*}{Theorem}
\theoremstyle{definition}

\newcommand{\bz}{\mathbb{Z}}

\newcommand{\bc}{\mathbb{C}} 

\newcommand{\coker}{\operatorname{coker}}
\newcommand{\fred}{S^2\times S^2}
\newcommand{\lan}{\left\langle}
\newcommand{\ran}{\right\rangle}

\newcommand{\rk}{\operatorname{rk}}
\newcommand{\fix}{\operatorname{Fix}}

\renewcommand{\hom}{\operatorname{Hom}}

\newcommand{\matrixb}{\left(\begin{smallmatrix} 0&1\\
	1&0\end{smallmatrix}\right)}

\title{Four-manifolds which admit $\bz_p\times\bz_p$ actions}
\author{Michael P. McCooey}
\address{Department of Mathematics and Statistics, McMaster University\\ 
1280 Main Street West\\
Hamilton, Ontario, Canada L8S 4K1}
\date{\today} 
\email{mmccooey@member.ams.org}

\subjclass{Primary 57S17, 57S25; Secondary 20J06}

\begin{document}
\begin{abstract}
We show that the simply-connected four-manifolds which admit locally
linear, homologically trivial $\bz_p\times\bz_p$ actions are
homeomorphic to connected sums of $\pm\bc P^2$ and $\fred$
(with one exception: pseudofree $\bz_3\times\bz_3$ actions on the
Chern manifold), and
also establish an equivariant decomposition theorem. 

This generalizes results from a 1970 paper by Orlik and Raymond about torus
actions, and complements more recent work of Fintushel, Yoshida, and Huck on
$S^1$ actions. In each case, the simply-connected four-manifolds which
support such actions are essentially the same.
\end{abstract}
\maketitle
\section{Introduction}
In 1970, Orlik and Raymond~\cite{OR} proved that any closed,
simply connected four-manifold which
admits a smooth, effective $S^1\times S^1$ action can be expressed as a
connected sum of copies of $\fred$, $\bc P^2$, and $-\bc P^2$. Later,
Fintushel~\cite{Fintushel77}  and Yoshida~\cite{Yoshida} each showed 
that the same conclusion holds for smooth $S^1$ actions.  In 1995, 
Huck~\cite{Huck} generalized this result to show that
the intersection form of a closed cohomology four-manifold $M$ with 
$H_1(M)=0$ on which $S^1$ acts must split as a sum of rank 
$1$ and $2$ forms $(\pm 1)$ and $\matrixb$, provided a certain regularity
condition holds near the fixed-point set of the action. Huck and 
Puppe~\cite{HuckPuppe} subsequently generalized further by removing the
restriction on $H_1(M)$. 

	Stated simply, Huck's approach is to study the equivariant cohomology 
of the singular set of an $S^1$ action using earlier techniques of 
Puppe~\cite{Puppe}, and thereby derive a characterization of the possible 
intersection forms. Related methods were used independently by the 
author~\cite{symmetrygroups} to study actions of finite nonabelian groups
on four-manifolds. Our methods actually simplify somewhat when the groups 
are abelian, and we apply them here to prove:
If $M$ is a closed four-manifold with $H_1(M) = 0$ which admits a locally
linear, homologically trivial action by $\bz_p\times\bz_p$ (with $p$ prime),
then the intersection form of $M$ splits as a sum of copies
of $(\pm 1)$ and $\matrixb$. 

	If this $M$ is simply connected, then by the work of 
Freedman~\cite{Freedman}, it must be homeomorphic to a connected sum of
copies of $\fred$, $\pm \bc P^2$, and perhaps a copy of 
$\pm\widehat{\bc P^2}$, where $\widehat{\bc P^2}$ denotes the manifold homotopy
equivalent to $\bc P^2$, but with non-vanishing Kirby-Siebenmann invariant.
We generalize an observation of Wilczy\'{n}ski~\cite{Wil2} to show (with
exactly one exception) that 
$KS(M)$ must vanish. As a corollary, we obtain an analogue of
Orlik and Raymond's result for $\bz_p\times\bz_p$ actions.

Finally, we discuss  the question of classifying the
actions themselves.  A complete classification would be very difficult, but
by combining our results with those of Orlik and Raymond, we prove an
equivariant version of the decomposition theorem which reduces
the general question to that of classifying the possible actions on $S^4$.

\section{The singular set of a $\bz_p\times\bz_p$ action}\label{ECcalcs}
Suppose $M$ is a closed, connected
 four-manifold with $H_1(M)=0$ and $b_2(M)\ge 1$.
If $G=\bz_p\times\bz_p $ acts on $M$, the set 
$\Sigma=\{x\in M\ |\ G_x\ne \{0\}\}$ is called 
the \emph{singular set} of the action.
We assume throughout the paper
 that the action is effective, locally linear, and
homologically trivial. 
By results of Edmonds~\cite{Aspects}, each $g\ne 0$ in $G$ has a fixed point
set consisting of isolated points and $2$-spheres, and each $2$-sphere
represents a nontrivial homology class. Our first task is to 
understand how the fixed-point sets of the cyclic subgroups $\bz_p\subset G$
fit together to form the overall structure of $\Sigma$.

Recall that the Borel
equivariant cohomology $H_G(X)$ of a $G$-space $X$ is the ordinary 
cohomology of the balanced product $EG\times_G X$. This balanced product
has a natural fibration over $EG/G=BG$, and the Leray-Serre spectral sequence
of the fibration is called the \emph{Borel spectral sequence}. 

The next lemma and the proposition which follows it appeared in slightly
different form in Edmonds~\cite{Aspects} and~\cite{Actions}, 
but we re-state them here for convenience:

\begin{lem}\label{sscollapsesforM} Suppose $G$ acts
homologically trivially on a closed four-manifold
$M$ with $H_1(M)=0$. If either
\begin{enumerate}

\item{ $H^2(M)$ contains a class  $u$ whose square generates
$H^4(M)$, and $H^3(G)$ has no $3$-torsion, or}
\item{$b_2(M)\ge 3$,} 
\end{enumerate}
then the Borel
spectral sequence $E(M)$ collapses with coefficients in $\bz$ or any field.
It follows that $H^*_G(M)$ is a free $H^*(G)$ module on $b_2(M)+2$ generators
corresponding to generators for $H^*(M)$.
\end{lem}
\begin{proof}
If $u\in H^2(M)$
has nonzero square, then, since $u^3=0$, $0=d_3(u^3)=3d_3(u)u^2$. But 
$E_3^{*, 4}$ is a free $H^*(G)$ module generated by $u^2$. So if 
$H^3(G)$ has no $3$-torsion, then $d_3(u)$ must be $0$. And then, of course,
$d_3(u^2)=0$, as well. Thus $E_2(M)=E_3(M)$. Since $d_5(u^2)=2ud_5(u)=0$,
the sequence collapses. 

Now suppose $b_2(M)\ge 3$. Then for each generator $u\in H^2(M)$,
there is a $v\in H^2(M)$ which is linearly independent of $u$ in $H^2(M)$,
and such that $uv=0$. Since the action of $G$ is homologically trivial, 
$E_2(M)$ is a free $H^*(G)$-module on generators corresponding  to 
those of $H^*(M)$, so in fact $u$ and $v$ must be independent in
$E_2(M)$, as well. But $d_3(uv)=ud_3(v)+vd_3(u)=0$. This is only possible
if $d_3(u)=d_3(v)=0$. $H^4(M)$ is generated by products of two-dimensional
classes, so $d_3^{*, 4}=0$, as well. It follows that $E_2(M)=E_3(M)=E_4(M)$.
The same argument shows that $E_4(M)=E_5(M)=E_{\infty}(M)$.

The conclusion about $H^*_G(M)$ follows immediately from 
tom Dieck~\cite[III.1.18]{tomDieck}.
\end{proof}

Whenever $G$ is a finite 
group and $S$ is a multiplicative, central subset
of the cohomology ring $H^*(G)$, we can define the ``$S$-singular set'' 
$\Sigma_S =\{x\in X| S\cap \ker r^*_{G_x}=\emptyset\}$. The fundamental
\emph{Localization Theorem} (See Hsiang~\cite{Hsiang} or 
tom Dieck~\cite{tomDieck}) then 
states that the localized restriction map 
$S^{-1}H^*_G(M)\rightarrow S^{-1}H^*_G(\Sigma_S)$ is an isomorphism.
Applying the Localization Theorem in specific cases requires careful
choice of $S$, based on knowledge of the restriction maps
from the cohomology of $G$ to that of its subgroups. But it can yield
useful information about the structure of $\Sigma$. 
We apply it to prove:

\begin{prop}\label{abelianfixedpointset}
Let $M$ be a closed four-manifold such that  $H_1(M)=0$,  and
suppose that either $b_2(M)\ge 3$, or $b_2(M)=2$ but the 
intersection form of $M$ is diagonalizable over $\bz$. 
If $G$ is a finite abelian group which acts locally linearly and 
homologically trivially on $M$, then the rank
of $G$ is at most 2, and $G$ has nonempty fixed-point set.
If $G=\bz_p\times\bz_p$, then $\fix(G)$ consists of exactly
$b_2(M)+2$ points.
\end{prop}

\begin{proof}
Suppose  first that $G=\bz_2\times\bz_2\times\bz_2$
acts on $M$. By Lemma~\ref{sscollapsesforM},
$H^*_G(M; \bz_2)$ is a free $H^*(G, \bz_2)$ module on $b_2(M)+2$ 
generators. Recall that $H^*(G; \bz_2)\cong \bz_2[a]\otimes\bz_2[b]
\otimes\bz_2[c]$, where $a,b,$ and $c$ generate $\hom(G,\bz_2)$. 
Since $H^*(G; \bz_2)$ is a polynomial ring, it contains no
zero-divisors, so it makes sense to localize at the set 
$S$ consisting of all of the nonzero elements. We check easily
that $S^{-1}H^*(G,\bz_2)\cong \bz_2$.

Now, each proper subgroup $H\subset G$ is the kernel of some nonzero
homomorphism $\varphi_H:G\rightarrow \bz_2$, and this $\varphi_H$,
viewed as an element of $H^1(G; \bz_2)$, restricts trivially
to $H^1(H; \bz_2)$. So the $S$-singular set $\Sigma_S$ contains only those
points fixed by all of $G$. By the Localization Theorem,
$\Sigma_S$ is nonempty.  Consideration of the isotropy representation
of $G$ at a fixed point $x_0$ shows that there must be $g, h\in G$ such that
$g$ fixes a two-dimensional subspace $V\subset T_{x_0}$, 
while $h|_V$ acts by 
$\left(\begin{smallmatrix}-1&0\\0&1\end{smallmatrix}\right)$.
But $V$ forms part of a $2$-sphere $S$ fixed by $g$. If $h$ reverses 
orientation on $V$, it also acts by $-1$ on $[S]\in H_2(M)$, contradicting
homological triviality.

If $G=\bz_2\times\bz_2$, the same argument shows that $\fix(G)$ 
contains $b_2+2$ points, but of course no contradiction ensues
from the isotropy representation of $G$.

If $p$ is odd, a similar argument applies, except in the case
where $b_2(M)=2$ and $p=3$. To ensure that $S$ 
is central in $H^*(G; \bz_p)$, we replace the
one-dimensional generators of $\hom(G, \bz_p)$,
with  their two-dimensional images under the Bockstein map, 
which generate the polynomial part of $H^*(G, \bz_p)$.

Finally, suppose $b_2(M)=2$ and $G=\bz_3\times\bz_3\times\bz_3$, with
generators $g, h$, and $k$. The Lefschetz fixed-point theorem implies
that $\chi(\fix(g))=4$. Thus $\fix(g)$ either 
contains at least one $2$-sphere, or consists of exactly four isolated 
points. In the first case, $G/\lan g \ran$ acts effectively on the 
sphere, which is impossible. In the second case, the action of $h$ on 
$\fix(g)$ must fix at least one point $x_0$. But $\lan g, h\ran$ cannot
act freely on the linking sphere to $x_0$, so some other element $g'$ fixes a 
$2$-sphere, and the argument proceeds as before.

\end{proof}

There are indeed actions of $\bz_2\times\bz_2$ 
on $\fred$, and pseudofree actions of 
$\bz_3\times\bz_3$ on $\bc P^2$ and $\widehat{\bc P^2}$, 
whose fixed point set is empty.  
Inspection shows that such actions are the 
only exceptions to the rule that $\fix(G)\ne \emptyset$,
and that if $\fix(G)\ne\emptyset$, 
then in fact it contains $b_2+2$ points.
Because our desired conclusion about the intersection form holds 
in the exceptional cases, we assume as a convenience in this 
section and the next that the fixed-point set is non-empty.
(Later, in the geometrical analysis, the assumption will be more essential,
and we will not take it for granted.)

Since the isotropy representation of $G$
at any fixed point splits as a sum of rank two real representations, 
each fixed point is included in exactly
two singular $2$-spheres. Since $G$ is abelian, $G$ acts on $\fix(g)$ for each
$g\in G$, so each sphere has a rotation action with fixed points at its north
and south poles. Thus $\Sigma$ contains a total of $b_2+2$ spheres $S_1, 
\ldots,S_{b_2+2}$, and
each path component of $\Sigma$ is a chain of such spheres arranged in a closed
loop. Since the action of $G$ on $\Sigma$ is just a rotation on each sphere,
$G$ acts trivially on $H^*(\Sigma)$.

\begin{lem}
If $p=2$, each $[S_i]$ represents a primitive homology class in $H_2(M, \bz)$.
\end{lem}

\begin{proof}
If each component of $\Sigma$ contains at least three spheres, then
each sphere intersects its neighbor geometrically once, and the 
claim follows. If some component contains exactly two spheres,
then each intersects the other twice. One of them might, a priori, represent 
a multiple of two in $H_2(M, \bz)$. 
But the theorem of Edmonds cited above implies
that it must be nontrivial in $H_2(M;\bz_2)$.
\end{proof}

If $p$ is odd, this argument does not suffice to  rule out certain
$[S_i]$ being multiples of 2. However, if $p$ is odd, then $2$-torsion
will not affect the cohomology calculations of the next section. 
The calculations of that section will show that $\Sigma$ is connected,
and then it will follow that each $S_i$ does, in fact, intersect its 
neighbor only once.

Our next goal is to show that the inclusion $H_2(\Sigma)\rightarrow
H_2(M)$ is (split) surjective. When we have shown this, it will follow that the
intersection form of $M$ is represented by the geometrical intersections
of the spheres in $\Sigma$.

From the cohomology long exact sequence of the pair $(M, \Sigma)$, we
extract:
$$0\rightarrow H^1(\Sigma)\rightarrow H^2(M, \Sigma)\rightarrow H^2(M)
\rightarrow H^2(\Sigma)\rightarrow H^3(M, \Sigma)\rightarrow 0.$$
A short diagram chase
shows that $G$ acts trivially on the relative cohomology groups. Let $N$
denote the number of path components of $\Sigma$, and $L$, the (integral) 
rank of $\coker H^1(\Sigma)\rightarrow H^2(M, \Sigma)$. 
As we have noted,
each $S_i$ represents an ``almost primitive'' homology class in $M$. 
More precisely: $H^3(M, \Sigma)\cong\bz^{L+2}\oplus T$, where $T=0$ if
$p=2$, and $2T=0$  if $p$ is odd.

We shall prove:
\begin{lem}
$L=0$,
\end{lem}
From which the claim about $H_2(\Sigma)\rightarrow H_2(M)$ is immediate.

\begin{proof}
Recall that 
$$H^*(\bz_2\times\bz_2; \bz)\cong\frac{\bz[\alpha_2, \beta_2]\otimes P[\mu_3]}
{\lan 2\alpha = 2\beta =2\mu =0, \mu^2=\alpha\beta^2+\alpha^2\beta\ran},$$
while for $p$ odd,
$$H^*(\bz_p\times\bz_p; \bz)\cong\frac{\bz[\alpha_2, \beta_2]\otimes E[\mu_3]}
{\lan p\alpha = p\beta =p\mu =0\ran}.$$

Let $\pi$ denote the projection $M\rightarrow M/G = M^*$.
The Borel spectral sequence of the pair $(M, \Sigma)$ has 
$$E_2^{i,j}(M, \Sigma)=H^i(G; H^j(M, \Sigma))\Rightarrow H^*_G(M, \Sigma).$$
On the other hand, $M - \Sigma$ is a free $G$-space, so 
$H_G^*(M, \Sigma)$ is canonically isomorphic to $H^*(M^*, \Sigma^*)$.

Since $(M^*, \Sigma^*)$ is
a relative manifold pair, Poincar\'{e} duality gives a commutative diagram:
\begin{equation*}\begin{CD}
@. H^3(M, \Sigma) @>{\cong}>> H_1(M - \Sigma)\\
@. @AA{\pi^*}A                                @VV{\pi_*}V\\
H_G^3(M, \Sigma)@>{\cong}>>H^3(M^*, \Sigma^*)      @>{\cong}>> H_1(M^* - \Sigma^*).
\end{CD}\end{equation*}

But $H_1(M - \Sigma)$ is generated by meridians to the spheres in $\Sigma$,
and each of these is a $p$-fold cover of its image in $H_1(M^* - \Sigma^*)$.
Thus $\pi_*$ is multiplication by $p$. Since the left-hand edge homomorphism
$E_2^{0,j}\rightarrow E_{\infty}^{0,j}$ of the Borel spectral
sequence is induced by the fiber inclusion
$j:(M, \Sigma)\rightarrow (M_G, \Sigma_G)$,  we can conclude that
$\coker(E_2^{0,3}\rightarrow E_{\infty}^{0,3})$ has exponent $p$.
In other words, no $\bz$ summand of $E_2^{0,3}$ supports more than one
non-zero differential. In rank counting arguments we can therefore treat
its integral rank as though it were a $\bz_p$ rank.

Notice also that, since $H^4(M^*, \Sigma^*)\cong \bz$, each nonzero
class in $E_2^{i,j}$ with $i+j\ge 4, i\ne 0,$ must be mortal. 

Consider the  terms of $E_2(M,\Sigma)$ indicated in Table~\ref{L=0table}: 

\begin{table}[h]
\caption{$E_2(M, \Sigma)$}\label{L=0table}
$\begin{array}{l|llllll}
H^4(M, \Sigma) & \bz & 0\\
H^3(M, \Sigma) & \bz^{2+L}\oplus T&0 \\
H^2(M, \Sigma) & \bz^{N+L}&0 &\bz_p^{2N+2L} \\
H^1(M, \Sigma) & \bz^{N-1}&0 & &\bz_p^{N-1} & \bz_p^{3N-3}\\
0&0&0&0&0&0&0\\
\hline
 & \bz & 0 & \lan \alpha, \beta\ran & \lan \mu\ran & \lan\alpha^2, \alpha\beta,
\beta^2\ran & \lan \mu\alpha, \mu\beta\ran  

\end{array}$\end{table} 

Now, elements of $E_2^{3,1}$ can only be killed by $d_2^{0,2}$, 
while $E_2^{2,2}$ can be killed 
either by $d_2^{0,3}$ or $d_2^{2,2}$. By the above observations, we have:
\begin{align*}
\rk E^{0,3} + \rk E^{4,1} &\ge \rk E^{2,2} + \rk E^{3,1}\\
(2+L) +(3N-3) &\ge (2N+2L) + (N-1),
\end{align*}
so $L=0$, as claimed. 
\end{proof}

We obtain a corollary which is dual, in some sense, to Edmonds's 
theorem~\cite[2.5]{Aspects}:
\begin{cor}
Suppose $G=\bz_p\times\bz_p$ acts as we have been assuming. Then the 
cohomology restriction map $H^2(M)\rightarrow H^2(\Sigma)$ is injective,
so $\Sigma$ represents all of $H_2(M)$.
\end{cor}

To better understand the geometry of the singular set, we also prove:

\begin{lem}
$N=1$, so $\Sigma$ is connected.
\end{lem}

\begin{proof}
From the homology spectral sequence of the covering $\pi:X - \Sigma
\rightarrow X^* - \Sigma^*$, we obtain a short exact sequence
$0\rightarrow \bz\times\bz\times T\stackrel{\pi_*}{\rightarrow} H_1(M^* - \Sigma^*)
\rightarrow \bz_p\times\bz_p\rightarrow 0$. As we have already seen,
$\pi_*$ is multiplication by $p$ in each factor. It follows that
$H_1(M^* - \Sigma^*)\cong H^3(M^*, \Sigma^*)$ is $p$-torsion-free.
So in fact all classes in $E_2^{i,j}$ with $i>0$ are mortal. In 
particular, $E_2^{1,2}$ must vanish, so $N\ge 2(N-2)$, so $N\le 2$.
Suppose for a contradiction that $N=2$. The Borel spectral
sequence then takes the following
form (each entry with $i>0$ is $\bz_p^k$ for some $k$, so to save space, 
we simply indicate its rank):

\begin{table}[h]
\caption{$E_2(M, \Sigma)$}
$\begin{array}{l|lccccc}
H^4(M, \Sigma) & \bz & 0& 2 & 1& 3& 2 \\
H^3(M, \Sigma) & \bz^{2}\oplus T&0& 4 & 2 & 6 & 4 \\
H^2(M, \Sigma) & \bz^{2}&0 &4 & 2 & 6 & 4 \\
H^1(M, \Sigma) & \bz    &0 &2 & 1 & 3 & 2 \\
0&0&0&0&0&0&0\\
\hline
 & \bz & 0 & \lan \alpha, \beta\ran & \lan \mu\ran & \lan\alpha^2, \alpha\beta,
\beta^2\ran & \lan \mu\alpha, \mu\beta\ran  

\end{array}$\end{table} 

There are generators $a\in H^1(M, \Sigma)$ and $b, c\in H^2(M,\Sigma)$ such
that $d_2(b)=\alpha a$ and $d_2(c)\ = \beta a$. By the multiplicative 
properties of the spectral sequence, this kills the entire row $j=1$,
\emph{except} $H^3(G; H^1(M, \Sigma)\cong\bz_p\cong\lan\mu\ran$. 

Now, $\ker d^{2,2}_2=\lan \beta b -\alpha c\ran$, so there is some $e\in
H^3(M, \Sigma)$ such that $d_2(e)=\beta b -\alpha c$. Since $E_3^{3,1}
=\lan \mu a \ran$ must also perish, there is $f\in H^3(M, \Sigma)$, 
independent of $e$, so that $d_3(f)=\mu a$. But then 
$d_3(\alpha f)=d_3(\beta f)
=0$, since $\mu\alpha a$ and $\mu\beta b$ were already killed by $d_2$. 
Now $\ker d_2^{2,3}$ has rank $2$, and $d_3^{2,3}=0$. But $d_2^{0,4}$ has
rank $\le 1$, so $E_{\infty}^{2,3}$ must have rank $\ge 1$. This is a 
contradiction, so $N=1$.
\end{proof}

Now that we know this, each $S_i$ definitely intersects each neighbor 
only once, so $T=0$ for odd $p$.

To summarize, we have shown:
\begin{prop}
Suppose $M$ is a closed,
topological four-manifold with $b_2(M)\ge 1$ and $H_1(M) = 0$, equipped with
an effective, homologically trivial, locally linear $\bz_p\times\bz_p$
action. With the exception of fixed-point free actions
which exist in the two cases,
\begin{enumerate}
\item{$b_2(M)=1$, $p=3$, and the action is pseudofree, or}
\item{$M$ has intersection form $\matrixb$ when $p=2$,}
\end{enumerate}
the singular set $\Sigma$ consists
of $b_2(M)+2$ spheres equipped with rotation actions, intersecting
pairwise at their poles, and arranged into a single closed loop. 
Each sphere represents a primitive class in $H_2(M; \bz)$, and together
these classes generate $H_2(M)$.
\end{prop}

\section{The intersection form}\label{intform}


Let $\sigma_1, \ldots, \sigma_{b_2+2}$ denote the fundamental classes
of $S_1, \ldots, S_{b_2+2}\in \Sigma$. As generators of $H_2(M)$, two can be
regarded as ``redundant''. If we eliminate one, 
we cut the loop of $\Sigma$. By removing another, 
we either disconnect  or shorten the remaining chain. Renumber the
remaining spheres, if necessary, as $S_1, \ldots, S_{b_2}$, and
call the result $\Sigma'$. Let $e_i=\sigma_i\cdot\sigma_i$. 
The matrix of intersections of the spheres, and therefore the intersection
form of $M$, as well, takes the form of one, or a sum of two, pieces
of the form
$$\begin{pmatrix}
e_1 & 1 & & & & & \\
1 & e_2 & 1 & & & & \\
& 1 & e_3 & 1 & & & \\
& &1 & \ddots & & & \\
& & & & & 1 & \\
& & & &1 & e_{k-1} & 1\\
& & & & & 1 & e_k 
\end{pmatrix}$$

Huck and Yoshida have already proven exactly the lemma we need about such 
matrices (See Huck~\cite[lemma 4.2]{Huck}): Each is equivalent to a sum 
of rank $1$ and $2$ pieces. Thus:
\begin{thm}\label{intformthm}
Suppose $M$ is a closed topological four-manifold with $H_1(M)=0$. 
Let $p$ be prime, and suppose $\bz_p\times\bz_p$ acts effectively,
locally linearly, and homologically trivially on $M$. Then the intersection
form of $M$ is a sum of copies of $(\pm 1)$ and $\matrixb$.
\end{thm}

\section{Vanishing of $KS(M)$}\label{KSvanishes}

Edmonds~\cite{Construction} showed that when $p$ is a prime greater 
than 3, locally linear, homologically trivial  $\bz_p$ actions exist on
every simply-connected four-manifold. The actions he constructs are 
\emph{pseudofree} -- i.e. with only isolated fixed points. In certain
cases, this is a necessary restriction. For example,  Wilczy\'{n}ski~\cite{Wil2} 
shows that if a homotopy $\bc P^2$ admits a $\bz_p$ action which fixes a 
two-sphere, then the two-sphere can be used to split off $\bc P^2$ as a 
connected summand of $M$, and it follows that $M$ is homeomorphic to $\bc P^2$.
In other words, $\widehat{\bc P^2}$ admits \emph{only} pseudofree actions.

For the remainder of the paper, we assume $M$ is simply connected.
We generalize  Wilczy\'{n}ski's construction to prove the following
corollary of Theorem~\ref{intformthm}:
\begin{thm} If $M$ is a closed, simply connected four-manifold which admits
a locally linear, homologically trivial $\bz_p\times\bz_p$ action,
then $M$ is homeomorphic to a connected sum of copies of $\pm \bc P^2$ and
$\fred$, or if $p=3$ and the action is pseudofree, perhaps to 
a single copy of $\pm\widehat{\bc P^2}$.
\end{thm} 

(In Theorem~\ref{equivariantversion}, we will establish a sharper result.)
\begin{proof}
For convenience, we continue to assume the action is not pseudofree,
or fixed-point-free in the case of $\fred$.
As above, let $\Sigma'$ denote the singular set with two homologically 
redundant spheres removed. 
For simplicity of notation, we assume $\Sigma'$
is connected; if it isn't, our argument will carry through on each piece.

It follows from the work of Freedman and Quinn~\cite[9.3]{FQ} (see 
also~\cite[1.2]{HLM}) that each 
$S^2\subset \Sigma$ has an equivariant normal bundle.
Thus $\Sigma'$ has a regular
neighborhood $N(\Sigma')$ which is  homeomorphic to the 
manifold obtained by plumbing together disk bundles $E(e_i)$, $i=1, 
\ldots, b_2$, over $S^2$ according to the graph 
$$A_{b_2} = e_1{\rule[2 pt]{2 em}{1 pt}}\ e_2{\rule[2 pt]{2 em}{1 pt}} 
\cdots{\rule[2 pt]{2 em}{1 pt}}\ \ e_{b_2-1}
{\rule[2 pt]{2 em}{1 pt}}\ e_{b_2}.$$
The boundary of such a plumbed manifold is a lens space $L$, and $|H_1(L)|$
is given by the determinant of the intersection matrix. In our case, the
matrix is unimodular, so $L$ is in fact a three-sphere. Thus $M'=N(\Sigma ')
\cup_{S^3}D^4$ is homeomorphic to a connected sum of copies of $\pm\bc P^2$ and
$\fred$. But $M'$ is also a connected summand of $M$ which carries all of its
homology. By Freedman and Quinn~\cite[10.3]{FQ}, $M'\cong M$.
\end{proof}

\section{On classifying $\bz_p\times\bz_p$ actions}

The argument of Orlik and Raymond on the classification of torus actions, 
specialized to the simply-connected case, can be summarized as follows:
The quotient space $M/T$ is a surface with boundary, and since $H_1(M)=0$,
it must be a disk. The boundary of $D$ consists of fixed points and arcs; 
the arcs can be labeled according to the corresponding isotropy subgroups
of $T$. Each arc lifts to a singular $S^2$ and each interior point of the
disk represents a principal orbit. They show that the quotient map in fact
admits an essentially unique section; thus the singular data in the quotient
space determine $M$ up to equivariant diffeomorphism. A calculation involving
the particular isotropy groups then shows that the quotient space splits in a 
way which lifts to an equivariant connected sum decomposition of $M$.

Up to an automorphism of $T$, listing the one-dimensional isotropy groups is 
equivalent to listing the Euler classes and (signed) intersection numbers 
of the singular $2$-spheres. This information is also available for 
$\bz_p\times\bz_p$ actions. To what extent does it classify them? We will 
show:

\begin{prop}\label{extensionlemma}
Assume the action is not one of the exceptional fixed-point-free 
cases.
\begin{enumerate}
\item{Each $\bz_p\times\bz_p$ action extends to a torus action in a 
regular neighborhood $\nu(\Sigma)$ of $\Sigma$.}
\item{$\nu(\Sigma)$ is $T$-equivariantly diffeomorphic to the singular set
of some smooth $T$-action on  $M$, but the given $T$-action
need not extend over $M$}.
\end{enumerate}
\end{prop}

\begin{proof}
We begin with a slight variant of the plumbing construction of the 
previous section: Let $t=b_2(M)+2$, and label the spheres consecutively around 
$\Sigma$ as $S_1, \ldots, S_t$. 
Let $x_i$ denote the ``north pole'' of $S_i$.
Choose orientations for each of the $S_i$, and let $\sigma_i$ denote the 
corresponding fundamental class. Finally, choose an orientation for $M$ and
let $\epsilon_i = \sigma_{i-1}\cdot\sigma_i$ denote the sign of the 
intersection at $x_i$. (When we considered $\Sigma'$ earlier, we implicitly
chose orientations to make each $\epsilon_i=+1$; here, because the spheres
are arranged in a closed loop, this might not be possible.)

With these conventions, $\nu(\Sigma)$ is obtained by plumbing together 
$D^2$-bundles $\xi_i$ over $S^2$, each with Euler class $e_i$, according
to the orientations given by the $\epsilon_i$. The plumbing graph is 
a circle, which we parameterize as $[\frac{1}{2},t+\frac{1}{2}]$ with 
the endpoints identified.
$\partial \nu$ can be thought of as a torus fiber bundle 
over the plumbing graph
with a fiber-preserving, free $\bz_p\times\bz_p$ action. It is not a priori
a principal bundle, but if the $\bz_p\times\bz_p$ action
on the fibers extends to a torus action, it will become one. With appropriate
smoothing around the plumbing points, the torus action will extend over 
$\nu(\Sigma)$, establishing the first part of the proposition.

The $T$-bundle over $[\frac{1}{2},t+\frac{1}{2}]$ can be assembled 
by gluing copies of 
$T\times [i-\frac{1}{2}, i+\frac{1}{2}]$ via attaching maps 
$\gamma_i$ which incorporate the clutching 
functions for the $\xi_i$, the coordinate switches at each plumbing point, 
and the orientations $\epsilon_i$. (See figure~\ref{plumbinggraph}, which
is intended to invite comparison with the diagrams in~\cite{OR}.) 
\begin{figure}
\caption{Torus bundle over the plumbing graph}
\label{plumbinggraph}\begin{center}\smallskip 
\resizebox{.5\textwidth}{!}{\includegraphics{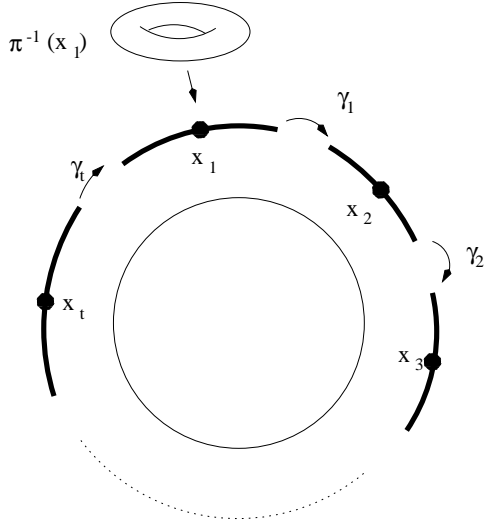}}
\end{center}\end{figure}
The maps are determined up to
isotopy by their $\pi_1(T)$ representations.
The clutching functions take the form $\left(\begin{smallmatrix}
1&0\\e_i&1\end{smallmatrix}\right)$; the coordinate switches are of
course $\left(\begin{smallmatrix}
0&1\\ 1&0\end{smallmatrix}\right)$, and the orientation changes,
$\left(\begin{smallmatrix}
1&0\\0&\epsilon_i\end{smallmatrix}\right)$. Together, such matrices
generate $GL(2,\bz)$, the structure group of the bundle.


The $\bz_p\times\bz_p$ action has well-defined rotation numbers in the fiber 
over $x_1$, so it extends to a torus action in that fiber. 
The gluing maps in $GL(2,\bz)$ define a trivialization of the bundle
over $T\times[\frac{1}{2},t+\frac{1}{2}]$ which is equivariant 
with respect to the $\bz_p\times
\bz_p$ action. Using them, the torus action extends along all of the fibers. 
The structure of the torus bundle is thus determined by the total 
gluing function
$\gamma:T\times \{t+\frac{1}{2}\}\rightarrow T\times \{\frac{1}{2}\}$;
with slight abuse of notation, we may write $\gamma = \gamma_t\circ\cdots
\circ\gamma_2\circ\gamma_1$. 

A compatibility condition is imposed by the 
existence of the $\bz_p\times\bz_p$ action -- namely, that the gluing
map $\gamma\in GL(2, \bz)$ must commute with order $p$ rotations
in each factor of $T\times 0$. We may analyze this requirement by lifting
to the universal cover $\pi:\widetilde{T}\rightarrow T$. 
A rotation $r$ lifts to a translation $\tau$. The requirement that
$\pi_*(\gamma^{-1}\tau\gamma)=\pi_*(\tau)=r$ means that the line 
spanned by each $\tau$ must be 
\begin{enumerate}
\item{
	Normalized by $\gamma$, if $p=2$. Since the total space of the 
bundle is the boundary of $\nu(\Sigma)$, it is orientable, so 
$\gamma$ is one of
$\pm\left(\begin{smallmatrix} 1 & 0\\0& 1\end{smallmatrix}\right)$ or
$\pm\left(\begin{smallmatrix}0 & 1\\- 1& 0\end{smallmatrix}\right),$}
\item{Centralized by $\gamma$, if $p>2$, which implies $\gamma =
\left(\begin{smallmatrix}1&0\\0&1\end{smallmatrix}\right)$.}
\end{enumerate}
 
In the latter case, $\gamma$ clearly commutes with the entire torus action,
so $\partial\nu$ supports the structure of a principal bundle. Even when $p=2$,
$\gamma$ must respect the base-fiber splitting of the bundle $\xi_t$ over
$S_t$, so $\pm\left(\begin{smallmatrix}0 & 1\\- 1& 0\end{smallmatrix}\right)$
is ruled out. We proceed to rule out $\gamma=\left(\begin{smallmatrix}-1 & 0\\
0& -1\end{smallmatrix}\right)$, also. If such a bundle is realized
on the singular set of a $\bz_2\times\bz_2$ action, let $\mu$ be a 
small meridional loop around $S_t$ in $M-\Sigma$. Then $\mu$ is homologous
to $-\mu$, and so $2\mu=0$ in $H_1(M-\Sigma)$. But $H_1(M-\Sigma)\cong
H^3(M, M-\Sigma)$ is torsion-free, as we saw in section~\ref{ECcalcs}.
It is generated by any pair of meridians to neighboring two-spheres in 
$\Sigma$.

This finishes the proof that $\partial\nu$ is a trivial principal $T$-bundle,
and hence also the proof of the first claim.  Our proof of the second claim
is constructive, based on Orlik and Raymond's model in the case of torus
actions.

The orbit space $A=\nu(\Sigma)/T$ is an annulus. Its outer boundary 
component $\partial_1A$ consists of $t$ fixed points separated by $t$ arcs 
whose stabilizers are copies of $S^1\subset T$. Its inner boundary $\partial_2A$
consists entirely of principal orbits. Adjoin a disk $D$ to $\partial_2A$. 
Because the torus bundle is trivial over $\partial_2A$, there is no obstruction
to lifting this adjunction to a $T$-equivariant gluing of $D^2\times T$
to $\partial\nu$. The resulting manifold, denoted $M'$, is simply connected
and has the same intersection form as $M$, so it is homeomorphic to 
$M$. 

Finally, an example of Hambleton, Lee, and Madsen (\cite{HLM}) shows that
$M$ and $M'$ need not be $\bz_p\times\bz_p$-equivariantly homeomorphic. 
They begin with a linear $\bz_p\times\bz_p$ action on $\bc P^2$, and
equivariantly connect sum a $\bz_p$-orbit of counterexamples to the
Smith conjecture in $S^4$ around one of the singular $2$-spheres. The
resulting space is still homeomorphic to $\bc P^2$, but the complement
of the singular set has nonabelian fundamental group. In the linear
example, $\bc P^2-\Sigma$ has the homotopy type of a torus.
\end{proof}

Let us call a $\bz_p\times\bz_p$ action
\emph{standard} if it is the restriction of a smooth torus action. It is fair
to say that the standard actions are completely understood. 
Proposition~\ref{extensionlemma}, together with the construction of
section~\ref{KSvanishes}, shows that we can equivariantly 
split off standard summands. If the two ``redundant'' two-spheres
are adjacent, then $M\cong M_{\text{standard}}\,\#\,S^4$, while if
there is no such choice of adjacent spheres, a two-step splitting
still yields $M\cong M_{\text{standard}}'\,\#\,M_{\text{standard}}''\,\#\,S^4$.
Because the standard actions extend to torus actions, Orlik and Raymond's
classification theorem applies to show that each splits further into 
``irreducible'' pieces.   This proves:

\begin{thm}\label{equivariantversion}
Let $M$ be a closed, simply-connected four-manifold with an effective,
locally linear, homologically trivial $\bz_p\times\bz_p$ action. 
Assume the action is not one of the fixed-point-free exceptions. Then $M$
admits an equivariant connected sum decomposition
$$M\cong S^4\#\, M_1 \#\ldots\#\, M_k,$$ where each $M_i$ is one of 
$S^4$, $\fred$, $\pm \bc P^2$, or $\bc P^2\# -\bc P^2$, equipped with a 
standard action. The action on the first $S^4$ summand need not be standard.
\end{thm}

Recall that the ``fixed-point-free exceptions'' are the pseudofree actions of
$\bz_3\times\bz_3$ on $\bc P^2$ and $\widehat{\bc P^2}$, and 
fixed-point-free actions of
$\bz_2\times\bz_2$ on $\fred$. Also note that Orlik and Raymond construct
examples of torus actions on $\bc P^2\# -\bc P^2$ which admit no
\emph{equivariant} connected sum decomposition.

As a consequence of this theorem, the general problem of 
classifying $\bz_p\times\bz_p$ actions on simply-connected
four-manifolds reduces to the  question of classifying actions on $S^4$. 
The latter is still, of course, very difficult.

\section{Questions}

Finally, we leave the reader with two questions:
\begin{enumerate}
\item{What can constructively be said about the classification of $\bz_p\times
\bz_p$ actions on $S^4$, in light of the possible knotting of 
the singular set?}
\item{Huck and Puppe~\cite{HuckPuppe} generalized Huck's earlier work 
on circle actions to the case $H_1(M)\ne 0$. Does Theorem~\ref{intformthm}
generalize similarly? It is worth noting that in the general case, the 
singular set need not contain spheres, as examples of free actions 
on $T^4$ easily show.}
\end{enumerate}

\bibliographystyle{abbrv}
\bibliography{mybiblio}

\end{document}